\documentclass[10pt]{article}
\usepackage[nohead,margin=1.0in]{geometry}
\usepackage{amssymb, amsmath, amsthm}
\usepackage{color}
\usepackage{graphicx}
\usepackage{dsfont}
\usepackage{epstopdf}
\graphicspath{{./pics/}}
\usepackage{booktabs}
\usepackage{threeparttable}
\usepackage{verbatim}
\newcommand{\bel}{\begin{equation} \label}
\newcommand{\ee}{\end{equation}}
\def\beq{\begin{equation}}
\def\eeq{\end{equation}}
\newcommand{\bea}{\begin{eqnarray}}
\newcommand{\eea}{\end{eqnarray}}
\newcommand{\beas}{\begin{eqnarray*}}
\newcommand{\eeas}{\end{eqnarray*}}

\newcommand{\re}{\mathfrak R}

\newcommand{\R}{\mathbb{R}}
\newcommand{\C}{\mathbb{C}} 
\newcommand{\N}{\mathbb{N}}

\renewcommand{\d}{{\rm d}}

\newtheorem{theorem}{Theorem}[section]
\newtheorem{lem}[theorem]{Lemma}

\newtheorem{proposition}[theorem]{Proposition}
\newtheorem{defn}[theorem]{Definition}
\newtheorem{remark}{Remark}[section]
\newtheorem{example}{Example}[section]
\numberwithin{equation}{section}

\renewcommand{\d}{\,\mathrm{d}}
\allowdisplaybreaks
\providecommand{\abs}[1]{\left\lvert#1\right\rvert}
\providecommand{\norm}[1]{\left\lVert#1\right\rVert}

\def\phi {\varphi}

\allowdisplaybreaks

\title{Inverse Coefficient Problem for One-Dimensional Subdiffusion\\ with Data on Disjoint Sets in Time}
\author{Siyu Cen\thanks{Department of Applied Mathematics, The Hong Kong Polytechnic University, Kowloon, Hong Kong, P.R. China (\texttt{siyu2021.cen@connect.polyu.hk}, \texttt{zhizhou@polyu.edu.hk})} \and Bangti Jin\thanks{Department of Mathematics, The Chinese University of Hong Kong, Shatin, New Territories, Hong Kong, P.R. China (\texttt{b.jin@cuhk.edu.hk, bangti.jin@gmail.com})}\and
Yavar Kian\thanks{Univ Rouen Normandie, CNRS, Normandie Univ, LMRS UMR 6085, F-76000 Rouen, France (\texttt{yavar.kian@univ-rouen.fr})}\and \'{E}ric Soccorsi\thanks{Aix-Marseille Universit\'e, Universit\'e de Toulon, CNRS, CPT, Marseille, France (\texttt{eric.soccorsi@univ-amu.fr})} \and Rachid Zarouf\thanks{Aix-Marseille Universit\'e, Laboratoire ADEF, Campus Universitaire de Saint-J\'er\^ome, 52 Avenue Escadrille Normandie Niemen, 13013 Marseille, France (\texttt{rachid.zarouf@univ-amu.fr})} \and Zhi Zhou\footnotemark[1]}

\date{}

\begin{document}
\maketitle

\begin{abstract}
In this work we investigate an inverse coefficient problem for the one-dimensional subdiffusion model, which involves a Caputo fractional derivative in time. The inverse problem is to determine two coefficients and multiple parameters (the order, and length of the interval) from one pair of lateral Cauchy data.
The lateral Cauchy data are given on disjoint sets in time with
a single excitation and the measurement is made on a time sequence located outside the support of the excitation. We prove two uniqueness results for different lateral Cauchy data. The analysis is based on the solution representation, analyticity of the observation and a refined version of inverse Sturm-Liouville theory due to Sini \cite{Sini:2004}. Our results heavily exploit the memory effect of fractional diffusion for the unique recovery of the coefficients in the model. Several numerical experiments are also presented to complement the analysis.
\end{abstract}

\section{Introduction}
In this work, we are concerned with an inverse coefficient problem for one-dimensional subdiffusion.
Fix the interval $(0,\ell)$, with $\ell>0$, $T>0$, let $q \in L^\infty(0,\ell)$ be a non-negative function and let $\rho\in L^\infty(0,\ell)$ satisfy the following condition
\begin{equation}\label{rho}
 0<c_0 \leq\rho(x) \leq C_0 <+\infty\quad\mbox{in } \Omega,
\end{equation}
with $c_0,C_0$ being two positive constants.
Also fix $g\in W^{1,1}(0,T)$, satisfying $g(0)=0$ and $\alpha \in (0, 1)$.
Now consider the weak solution  $u$ of the following initial boundary value problem (IBVP):
\begin{equation}\label{eq1}
\begin{cases}
\rho(x)\partial_t^{\alpha}u -\partial_x^2 u+q(x)u =  0, & \mbox{in }(0,\ell)\times(0,T),\\
  u(0,t)= g(t),\ u(\ell,t)= 0, & \mbox{on } (0,T), \\
u=0, & \mbox{in } (0,\ell)\times \{0\}.
\end{cases}
\end{equation}
In the model, the notation $\partial_t^\alpha u$ denotes the Caputo fractional derivative in time of order $\alpha\in(0,1)$, defined by (see, e.g., \cite[p. 92]{KilbasSrivastavaTrujillo:2006} or \cite[p. 41]{Jin:2021})
\begin{equation}
    \partial_t^\alpha u (t) = \frac{1}{\Gamma(1-\alpha)}\int_0^t(t-s)^{-\alpha}u'(s) \d s,
\end{equation}
where $\Gamma(z)=\int_0^\infty s^{z-1}e^{-s}\d z$ for $\Re(z)>0$ denotes Euler's Gamma function.
The model \eqref{eq1} arises in the mathematical modeling of subdiffusion processes, in which the mean squared particle displacement grows sublinearly with the time, as opposed to linear growth in normal diffusion. Microscopically, such transport processes can be described by continuous time random walk in which the waiting time between two consecutive jumps follows a heavy-tailed distribution (see, e.g., \cite{Meerschaert:2019} and \cite[Chapter 1]{Jin:2021}), and the probability density function of the particle appearing at location $x$ and time $t$ satisfies the model \eqref{eq1}. Subdiffusion has been observed in diverse applications, e.g., thermal diffusion in medium with fractal geometry \cite{Nigmatulin:1986}, protein transport within membranes \cite{Kou:2008}, ion dispersion in column experiments \cite{HatanoHatano:1998} and dispersion in heterogeneous
aquifer \cite{AdamsGelhar:1992}. See the surveys \cite{MetzlerJeon:2014,MetzlerKlafter:1998} for long lists of applications in physics and biology.

In Proposition \ref{pp1}, we prove  that problem \eqref{eq1} admits a unique weak solution $u\in L^1(0,T;H^{2}(0,\ell))$. Moreover, under suitable assumptions on the Dirichlet data $g$, we prove that there exists $\delta\in(0,T)$ such that the map $t\mapsto \partial_xu(0,t)$ belongs to  $C((T-\delta,T))$. By fixing $(t_n)_{n\in\mathbb N}$ an arbitrary increasing  sequence of $(T-\delta,T)$ that converges to $T$, we consider the  following inverse coefficient problem:   Determine simultaneously the coefficients $\rho$ and $q$ and the parameters $\alpha$ and $\ell$ from the knowledge of pointwise flux measurement $\partial_xu(0,t_n)$ for $n\in\mathbb N$. In addition, we also consider a variant of the inverse problem with the Neumann data at $x=\ell$.

Several variants of the inverse coefficient problem have been investigated \cite{JinZhou:2021ip,JingYamamoto:2023,RundellYamamoto:2018,RundellYamamoto:2023}. The unique recovery of the potential $q$ alone from the lateral flux measurement was studied by Rundell and Yamamoto \cite{RundellYamamoto:2018,RundellYamamoto:2023}. Later Jing and Yamamoto \cite{JingYamamoto:2023} extended the results to simultaneously recover the potential $q$, the order $\alpha$, the constants in the Robin boundary condition, and the initial value from the lateral Cauchy data at $x=0,\ell$, when the initial value satisfies a suitable nonvanishing condition. 
When the initial condition and (time-independent) source are unknown, Jin and Zhou \cite{JinZhou:2021ip} showed the unique recovery of the potential $q$ and the order $\alpha$ when the boundary excitation $g$ is judiciously chosen. All these works \cite{JinZhou:2021ip,JingYamamoto:2023,RundellYamamoto:2018,RundellYamamoto:2023} are concerned with the 1D model, and the key tool in the analysis is Gel'fand-Levitan theory for the inverse Sturm-Liouville problem \cite{PT}.  The unique determination of several coefficients was addressed in \cite{Kian:2022ip,KianLiLiuYamamoto:2021} for multi-dimensional problems from suitable choice of the excitation based on the works  \cite{AS,CY} (see also the work \cite{HLLZ} for a similar approach with internal measurement). Cen et al \cite{CenJinZhou:2023} developed a multi-dimensional extension of the results in \cite{JinZhou:2021ip}  for piecewise constant diffusion coefficient with one polygonal / circular inclusion, and employed the level set method  to recover the inclusion shape.  More broadly, the present work lies in inverse coefficient problems for subdiffusion, which have been extensively studied; See the reviews \cite{JinRundell:2015,LiYamamoto:2019coeff} for further results on unique recovery of coefficients from over-posed data. 

Most of the above mentioned results have treated inverse coefficient problems with an overlap between excitation and measurements. In contrast, we investigate the inverse problem with data on disjoint sets in time with excitation and measurement located on different sets in time. Namely, the excitation $g$ will be applied along the time interval $(0,T-\delta)$ while the measurement will be made on a sequence of time lying in $(T-\delta,T)$. This setting is not only useful for  applications but also very challenging, and so far only few mathematical results in that direction are available. We refer to Section \ref{sec:main} for the precise statements of the main results, and further discussions on relevant results. 
Moreover, we present numerical results to show the feasibility of the reconstruction using the Levenberg-Marquardt method \cite{Levenberg:1944,Marquardt:1963}.

The rest of the paper is organized as follows. In Section \ref{sec:main} we describe the main theoretical results and discuss them in connection with existing results. In Section \ref{sec:prelim} we give preliminary results on the model \eqref{eq1}, e.g., the solution representation and analyticity of the data. In Section \ref{sec:proof} we present the proofs of the main results. Finally in Section \ref{sec:numer}, we provide numerical results to complement the theoretical analysis.

\section{Main results and discussions}\label{sec:main}

In this section we present the main theoretical results of the work. The first result gives the unique recovery of two coefficients and two parameters for problem \eqref{eq1}. The proof of the theorem is given in Section \ref{sec:proof}. 
\begin{theorem}
\label{t1} 
For $j=1,2$, let $\alpha_j \in(0,1)$, $\ell_j \in (0,\infty)$, $q_j\in L^\infty(0,\ell_j)$ be non-negative, and $\rho_j\in L^\infty(0,\ell_j)$ be piecewise constant and fulfill condition \eqref{rho} with $\rho=\rho_j$. Let $g\in W^{1,1}(0,T)$ be not everywhere zero and satisfy $g(0)=0$ and
\bel{t1a}
\exists \delta \in (0,T),\ t \in (T-\delta,T) \quad \Longrightarrow\quad g(t)=0.
\ee
Denote by $u_j$, $j=1,2$, the solution to \eqref{eq1} with $(\alpha,\ell,q,\rho)=(\alpha_j, \ell_j,q_j, \rho_j)$.
Then the map $t\mapsto \partial_xu_j(0,t)\in  C((T-\delta,T))$ and, for any increasing sequence $(t_n)_n \in (T-\delta,T)^{\mathbb N}$ converging to $T$, we have
\bel{t1b}
\partial_xu_1(0,t_n)=\partial_xu_2(0,t_n),\ n\in\mathbb N\quad
\Longrightarrow \quad (\alpha_1,\ell_1,\rho_1,q_1)=(\alpha_2,\ell_2,\rho_2,q_2).
\ee
\end{theorem}

Next we study a variant of the concerned inverse problem. We pick $\alpha \in (0, 1)$ and $h\in W^{1,1}(0,T)$, satisfying $h(0)=0$, and we consider the following IBVP:
\begin{equation}\label{eq2}
\begin{cases}
\rho(x)\partial_t^{\alpha}w -\partial_x^2 w+q(x)w =  0, & \mbox{in }(0,\ell)\times(0,T),\\
w(0,t)= h(t),\ \partial_xw(\ell,t)= 0, & \mbox{on }(0,T), \\
w=0, & \mbox{in } (0,\ell)\times \{0\}.
\end{cases}
\end{equation}
We assume that there exists an open not-empty  subset $\omega\subset (0,\ell)$ and a constant $c>0$ such that
\begin{equation}\label{q}
q(x)\geq c,\quad x\in\omega.
\end{equation}

Then we can prove the following unique recovery result.
\begin{theorem}\label{t2} Let the condition of Theorem \ref{t1} be fulfilled with $\rho_j\in C^1([0,1])$ and $q_j$ satisfying  condition \eqref{q} with $q=q_j$, $j=1,2$.
Let also  one of the following conditions:
{\rm(i)} $\rho_1=\rho_2\equiv1$, and {\rm (ii)} $q_1=q_2$ and $\rho_1(0)=\rho_2(0)$, be fulfilled.
Let $h\in W^{1,1}(0,T)$ be not everywhere zero and satisfy $h(0)=0$ and
\bel{t2b}\exists \delta \in (0,T),\ t \in (T-\delta,T) \quad\Longrightarrow\quad h(t)=0.\ee
Denote by $w_j$, $j=1,2$, the solution to \eqref{eq2}, with  $(\alpha,\ell,q,\rho)=(\alpha_j, \ell_j,q_j, \rho_j)$.
Then the map $t\mapsto \partial_xw_j(0,t)\in  C((T-\delta,T))$ and, for any increasing sequence $(t_n)_n \in (T-\delta,T)^{\mathbb N}$ converging to $T$, we have 
\bel{t2b} \partial_xw_1(0,t_n)=\partial_xw_2(0,t_n),\ n\in\mathbb N\quad
\Longrightarrow \quad (\alpha_1,\ell_1,\rho_1,q_1)=(\alpha_2,\ell_2,\rho_2,q_2).\ee
\end{theorem}

One distinct feature of Theorems \ref{t1} and \ref{t2} lies in the fact that we have established the simultaneous determination of coefficients and other parameters from the observational data on disjoint subsets of the time interval $(0,T)$. Indeed, while the excitation $g$ (respectively $h$) is supported in $(0,T-\delta)$, the measurements are made at the sequence in time $(t_n)_n \in (T-\delta,T)^{\mathbb N}$ and there is no overlap between the excitation and the measurement. Such results are not only significant in terms of practical applications, since in many situations, it is often impossible to simultaneously apply excitation and make the measurement, but also in terms of mathematical contribution since only few results of determination of coefficients from data on disjoint sets are available. To the best of our knowledge, the only  other result proving determination of coefficients from single measurement and data on disjoint  set in time is given in 
\cite{Kian:2022ip} where the author combined the approach of \cite{KianLiLiuYamamoto:2021} for coefficients determination and the memory effect for subdiffusion \cite{KJ}. However, the result of  \cite{Kian:2022ip} requires not only the specific class of excitation (which is not easy to realize numerically), inspired by the prior works 
\cite{AS,CY}, but also extra boundary measurements (when compared with the work \cite{KianLiLiuYamamoto:2021}). In this work, we show for the first time that the memory effect of subdiffusion can be applied to single boundary measurements associated with a very general class of excitation and less measurements than \cite{JinZhou:2021ip} for the determination of coefficients in 1D subdiffusion. This confirms the observation for inverse source problems \cite{JK,Yamamoto:2023isp} where the source term was recovered from \textit{a posteriori} boundary measurement. In addition, we require the measurements only be a discrete sequence in time $(t_n)_n \in (T-\delta,T)^{\mathbb N}$, while all other similar results consider measurement on an interval in time.

Also  Theorems \ref{t1} and \ref{t2}  substantially improve the existing result \cite{JinZhou:2021ip} for the 1D problem by reducing the amount of required observational data and yet recovering more unknown coefficient / parameters, including the density $\rho$ and the length $\ell$ of the interval. This kind of simultaneous determination of different parameters can be compared  with \cite{Kian:2022ip,KianLiLiuYamamoto:2021} where similar results can be found for the multi-dimensional problems. However, in contrast to \cite{Kian:2022ip,KianLiLiuYamamoto:2021}, Theorems \ref{t1} and \ref{t2} are stated with a general class of excitations only subjected to the conditions \eqref{t1a} and \eqref{t2b}, whereas the excitations in \cite{Kian:2022ip,KianLiLiuYamamoto:2021} are far more specialized. Thus, our results are much more suitable than those of \cite{Kian:2022ip,KianLiLiuYamamoto:2021} for numerical computation and practical applications. Actually, in contrast to \cite{Kian:2022ip,KianLiLiuYamamoto:2021}, we present in Section \ref{sec:numer}  a numerical reconstruction method along with different examples based on the results of Theorems \ref{t1} and \ref{t2}.

\section{Analysis of the direct problem}\label{sec:prelim}

In this section we establish that each of the IBVPs \eqref{eq1} and \eqref{eq2} admits a unique weak solution enjoying time-analytic regularity. Since both problems can be studied in a similar fashion, we shall solely focus on the IBVP \eqref{eq1} and only state the corresponding result for problem \eqref{eq2} briefly.

\subsection{Preliminaries}
\label{sec-prem}

First we recall the definition of weak solutions  of problems \eqref{eq1} and \eqref{eq2}, following \cite{KY1}; see also \cite[Chapter 6]{Jin:2021}, \cite{KubicaYamamoto:2020} and  \cite{SY} for related discussions.
The notation $\langle\cdot,\cdot\rangle_{H^{-1}(0,\ell),H_0^1(0,\ell)}$ denotes duality pairing between the spaces $H^{-1}(0,\ell)$ and $H_0^1(0,\ell)$.
\begin{defn}
\label{d1} 
A function $u\in W^{1,1}(0,T;H^{-1}(0,\ell))\cap L^1(0,T;H^2(0,\ell))$ is said to be a weak solution of problem \eqref{eq1} if the following conditions are fulfilled:
{\rm(a)} for all $\psi\in H^1_0(0,\ell)$ and a.e. $t\in(0,T)$ we have
\bel{d1a}\left\langle \partial_t^\alpha u(\cdot,t),\psi\right\rangle_{H^{-1}(0,\ell),H^1_0(0,\ell)}+\int_0^\ell\rho^{-1}(x)(-\partial_x^2 u(x,t) +q(x) u(x,t))\psi(x)\d x=0;
\ee
{\rm(b)} for a.e. $t\in(0,T)$, $u(\ell,t)= 0$,  $u(0,t)= g(t)$;
and {\rm(c)} for all $\psi\in H^1_0(0,\ell)$, $\left\langle  u(\cdot,0),\psi\right\rangle_{H^{-1}(0,\ell),H^1_0(0,\ell)}=0$.
\end{defn}

\begin{defn}\label{d2}A function $w\in W^{1,1}(0,T;H^{-1}(0,\ell))\cap L^1(0,T;H^2(0,\ell))$ is said to be a weak solution of problem \eqref{eq2} if the  conditions {\rm (a)} and {\rm (c)} of Definition \ref{d1} are fulfilled  with $u=w$, and  for a.e. $t\in(0,T)$, $\partial_x w(\ell,t)= 0$,  $w(0,t)= h(t)$.\end{defn}

Now we recall useful notations. $L^2(0,\ell;\rho \d x)$ denotes 
the set  of measurable real valued functions $f$ such that
$\int_0^\ell |f(x)|^2\rho(x)\d x < \infty.$
Endowed with the scalar product 
$\left\langle f,g \right\rangle=\int_0^\ell f(x)g(x)\rho(x)\d x$, $L^2(0,\ell;\rho \d x)$
is a Hilbert space. Moreover,  condition \eqref{rho} ensures that $L^2(0,\ell;\rho \d x)=L^2(0,\ell)$. In the sequel, we write $\norm{f}$ for $\langle f , f \rangle^{1 \slash 2}$ for all $f \in L^2(0,\ell;\rho \d x)$.
Consider the unbounded linear operator in $L^2(0,\ell;\rho \d x)$ acting on its domain $D(A)=H^1_0(0,\ell)\cap H^2(0,\ell)$ defined as
$$Af=-\rho^{-1}(f''+qf),\quad f\in D(A).$$
Then  $A$ is selfadjoint in $L^2(0,\ell;\rho \d x)$ and its spectrum consists of an increasing sequence of simple positive eigenvalues 
$(\lambda_{n})_n$, see, e.g., \cite[Chap. 2, Theorem 2]{PT}. We pick an orthonormal basis $\{\phi_k:\ k\in\mathbb N\}$ in $L^2(0,\ell;\rho \d x)$ of real-valued eigenfunctions of $A$, such that
$$ A \varphi_k = \lambda_k \varphi_k,\quad k \in \N. $$
Evidently, 
$
\phi_{k}'(0) \neq 0,\ k\in\mathbb N,
$
since, otherwise, $\varphi_k$ would be zero everywhere \cite[Chap. 1, Corollary 1]{PT} which contradicts the fact that it is normalized. 
Next, we consider the unique $H^2(0,\ell)$-solution $v$ of the problem 
\begin{equation}\label{BVP}
\begin{cases}
- v''(x)+q(x)v(x) =  0, & x \in (0,\ell),\\
 v(\ell)= 0,\ v(0)= 1. & 
\end{cases}
\end{equation}

We have the following technical result.

\begin{lem}
\label{l1} 
Let $v$ be the unique $H^2(0,\ell)$-solution to \eqref{BVP}. Then we have
\begin{equation*}
\left\langle v,\phi_k\right\rangle=-\frac{\phi_{k}'(0)}{\lambda_k},\quad k\in\mathbb N,
\quad\mbox{and}\quad
\sum_{k=1}^\infty\frac{(\phi_{k}'(0))^2}{\lambda_k^2}<\infty.
\end{equation*}
\end{lem}
\begin{proof}
Multiplying by $\phi_k$ on both sides of the first equation in \eqref{BVP} and then integrating over $(0,\ell)$, gives
$$ -\int_0^\ell v''(x) \varphi_k (x) \d x + \int_0^\ell q(x) v(x) \varphi_k (x) \d x = 0. $$
Now, by integrating by parts in the first integral and using the identity
$-\varphi_k'' + q \varphi_k = \rho \lambda_k \varphi_k$ in 
$L^2(0,\ell)$, we get the desired identity. Finally, the estimate follows readily from this and the Bessel-Parseval theorem.
\end{proof}

\subsection{Well-posedness}
We shall prove with the aid of Lemma \ref{l1} that the problem \eqref{eq1} is well-posed, that is to say that it possesses a unique weak solution in the sense of Definition \ref{d1}. We will use frequently the two-parameter Mittag-Leffler function $E_{\alpha,\beta}$, $(\alpha,\beta) \in \R^2$ (see e.g. \cite[Section 3.1]{Jin:2021} or \cite[Eq. (1.56)]{P}), defined by
\begin{equation*}
E_{\alpha,\beta}(z) = \sum_{n=0}^\infty \frac{z^n}{\Gamma(n\alpha+\beta)},\quad z\in \mathbb{C}.
\end{equation*}

The first main result of this section is as follows.
\begin{proposition}
\label{pp1} 
Let $g\in W^{1,1}(0,T)$ satisfy $g(0)=0$, $q\in L^\infty(0,\ell)$ be non-negative, and let $\rho\in L^\infty(0,\ell)$ fulfill condition \eqref{rho}. Then, there exists a unique weak solution $u$ to problem \eqref{eq1} in the sense of Definition \ref{d1}, which is expressed by
\bel{pp1a}
u(x,t)=-\sum_{k=1}^\infty \phi_{k}'(0) \phi_{k}(x) \left(\int_0^t(t-s)^{\alpha-1}E_{\alpha,\alpha}(-\lambda_{k}(t-s)^{\alpha})g(s)\d s\right),\quad x\in (0,\ell),\  t\in(0,T).\ee
Moreover, for a.e. $t\in(0,T)$ we have
\bel{pp1b}
\partial_xu(0,t)=v'(0)g(t)+\sum_{k=1}^\infty \frac{(\phi_{k}'(0))^2}{\lambda_k} \left(\int_0^tE_{\alpha,1}(-\lambda_{k}(t-s)^{\alpha})g'(s)\d s\right),
\ee
where $v$ solves  problem \eqref{BVP} and the sequence on the right hand side converges in $L^1(0,T)$.
\end{proposition}
\begin{proof} 
Repeating the argument of \cite[Theorem 2.5.]{KY1} yields that the map $u$ given by \eqref{pp1a} is   the unique weak solution of \eqref{eq1} in the sense of Definition \ref{d1}.\footnote{Note that in \cite[Theorem 2.5.]{KY1} the authors considered \eqref{eq1} with the domain $\Omega$ being an open bounded space of $\R^d$ with $d\geq2$ but their argumentation can be easily extended to $\Omega=(0,\ell)$.} Thus, we only need to show that the representation \eqref{pp1b} holds. Integrating by parts in \eqref{pp1a} and applying \cite[Lemma 3.2]{SY} and Lemma \ref{l1} lead to
\begin{align}
u(\cdot,t) & = -\sum_{k=1}^\infty \frac{\phi_{k}'(0)}{\lambda_k} g(t) \phi_{k}+\sum_{k=1}^\infty \frac{\phi_{j,k}'(0)}{\lambda_k}\left(\int_0^tE_{\alpha,1}(-\lambda_{k}(t-s)^{\alpha})g'(s)\d s\right) \phi_{k} \nonumber \\
& =  g(t) v+\sum_{k=1}^\infty \frac{\phi_{j,k}'(0)}{\lambda_k} \left(\int_0^tE_{\alpha,1}(-\lambda_{k}(t-s)^{\alpha})g'(s)\d s\right) \phi_{k},\quad  t\in(0,T).
\label{pp1c}
\end{align}
As we have $v\in H^2(0,\ell)$, it suffices to show that the above series lies in $L^1(0,T;H^{2}(0,\ell))$. 
Since the embedding $D(A) \hookrightarrow H^{2}(0,\ell)$ is continuous, $L^1(0,T;D(A))$ is continuously embedded into $L^1(0,T;H^{2}(0,\ell))$. Therefore it is enough to prove that the sequence of general term
$$
f_N(x,t)=\sum_{k=1}^N \frac{\phi_{j,k}'(0)}{\lambda_k} \left(\int_0^tE_{\alpha,1}(-\lambda_{k}(t-s)^{\alpha})g'(s)\d s\right)\phi_{k},\quad x \in (0,\ell),\ t\in(0,T),\ N\in\mathbb N, 
$$
converges in $L^1(0,T;D(A))$. Moreover, since 
$L^1(0,T;D(A))$ is a Banach space, we can proceed by showing that $(f_N)$ is a Cauchy sequence in $L^1(0,T;D(A))$. Indeed, for any $M,N\in\mathbb{N}$ with $M< N$, we have 
\begin{align}  \norm{f_N-f_M}_{D(A)} \leq  \int_0^t  \norm{\sum_{k=M+1}^N \frac{\phi_{j,k}'(0)}{\lambda_k} E_{\alpha,1}(-\lambda_{k}(t-s)^{\alpha})\phi_{k}}_{D(A)} \abs{g^\prime(s)}\d s,\quad t \in (0,T). \label{e1}
\end{align}
Further, since for any $\lambda>0$ and $t>0$,
\bel{essess}
\abs{E_{\alpha,1}(-\lambda t^{\alpha})}\leq C\lambda^{-1}t^{-\alpha},\quad 
\ee
for some constant $C$ independent of $\lambda$ and $t$ \cite[Theorem 3.6]{Jin:2021}, we have
\begin{align*}
\norm{\sum_{k=M+1}^N \frac{\phi_{k}'(0)}{\lambda_k} E_{\alpha,1}(-\lambda_{k}(t-s)^{\alpha})\phi_{k}}_{D(A)} 
& = \left( \sum_{k=M+1}^N \abs{\phi_{k}'(0)}^2  E_{\alpha,1}(-\lambda_{k}(t-s)^{\alpha})^2 \right)^{\frac{1}{2}} \\
& \leq C  \left(\sum_{k=M+1}^N \frac{|\phi_{k}'(0)|^2}{\lambda_k^2}\right)^{\frac{1}{2}} (t-s)^{-\alpha}.
\end{align*}
This, \eqref{e1} and
Young's convolution inequality (which can be applied here since  $\alpha \in (0,1)$) lead to 
\begin{align*}
\norm{f_N-f_M}_{L^1(0,T;D(A))}
\leq& C \left(\sum_{k=M+1}^N \frac{|\phi_{k}'(0)|^2}{\lambda_k^2}\right)^{\frac{1}{2}}\norm{ \int_0^t (t-s)^{-\alpha}|g'(s)|\d s}_{L^1(0,T)}\\
\leq& C \left(\sum_{k=M+1}^N \frac{|\phi_{k}'(0)|^2}{\lambda_k^2}\right)^{\frac{1}{2}}\left(\int_0^T t^{-\alpha}\d t\right) \norm{g'}_{L^1(0,T)}.
\end{align*}
In view of Lemma \ref{l1}, this entails
$$\lim_{M,N\to+\infty}\norm{f_N-f_M}_{L^1(0,T;D(A))}=0,$$
which proves that $(f_N)$ converges in $L^1(0,T;D(A))$ and hence $u\in L^1(0,T;H^2(0,\ell))$. Moreover, since 
$(f_N)$ is convergent in $L^1(0,T;H^2(0,\ell))$ and since $H^2(0,\ell)$ is continuously embedded into $C^1([0,\ell])$, by Sobolev embedding theorem, the series on the right hand side of \eqref{pp1c} converges in $L^1(0,T;C^1([0,\ell]))$. Therefore, for a.e. $t\in(0,T)$, we have
\begin{align*}
   &\partial_x\left(\sum_{k=1}^\infty \frac{\phi_{k}'(0)}{\lambda_k}\left(\int_0^tE_{\alpha,1}(-\lambda_{k}(t-s)^{\alpha})g'(s)\d s\right)\phi_k(x)\right)\big|_{x=0}\\
   =&\sum_{k=1}^\infty \frac{(\phi_{k}'(0))^2}{\lambda_k}\left(\int_0^tE_{\alpha,1}(-\lambda_{k}(t-s)^{\alpha})g'(s)\d s\right),
\end{align*}
and this last sequence converges in $L^1(0,T)$. This proves that  \eqref{pp1b} holds true for a.e. $t\in(0,T)$.
\end{proof}

By Proposition \ref{pp1}, the Neumann trace $\partial_x u(0,t)$ at $x=0$ of the solution $u$ to \eqref{eq1} is well-defined for a.e. $t \in (0,T)$, 
and expressed by \eqref{pp1b}. Still, this is not enough to rigorously define the data $( \partial_x u(0,t_n) )_{n \in \N}$, where $(t_n)$ is 
an arbitrary sequence in $(0,T)$ converging to $T$, used in Theorem \ref{t1}. Nonetheless, this can be remedied by assuming that the Dirichlet data $g$ at $x=0$ in \eqref{eq1} vanishes in the vicinity of $T$. 

\subsection{Time-analyticity}

In this part we build a time-analytic extension of the solution $u$ to \eqref{eq1} under condition \eqref{t1a}.

\begin{proposition}
\label{pp2} 
Let $q$ and $\rho$ be the same as in Proposition \ref{pp1}. Let $g \in W^{1,1}(0,T)$ satisfy $g(0)=0$ and condition \eqref{t1a}. Then the mapping
 $(T-\delta,T) \ni t \ \mapsto \partial_x u(0,t)$, where $\delta \in (0,T)$ is defined by \eqref{t1a}, extends to an analytic function 
 \bel{pp2a}
t\mapsto \sum_{k=1}^\infty \frac{(\phi_{k}'(0))^2}{\lambda_k} \left(\int_0^{T-\delta}E_{\alpha,1}(-\lambda_{k}(t-s)^{\alpha})g'(s)\d s\right)\quad \mbox{in }(T-\delta,\infty).
\ee
\end{proposition}
\begin{proof} 
From \eqref{t1a} and \eqref{pp1b}, we have for a.e. $t\in(T-\delta,T)$,
\bel{e3}
\partial_xu(0,t) 
=\sum_{k=1}^\infty \frac{(\phi_{k}'(0))^2}{\lambda_k} \left( \int_0^{T-\delta}E_{\alpha,1}(-\lambda_{k}(t-s)^{\alpha})g'(s)\d s \right).
\ee
 Hence, it suffices to show that the series on the right-hand side of \eqref{e3} depends analytically on $t$ in $(T-\delta,\infty)$.
We proceed by extending the function defined by \eqref{pp2a} to a suitable cone of the complex plane $\mathbb{C}$.  
By \cite[Theorem 3.2]{Jin:2021} or \cite[Theorem 1.6]{P},  for all $\mu \in (\frac{\pi \alpha}{2}, \pi \alpha)$, there exists $C>0$ such that
\bel{e2} 
\forall z \in \C,\ \mu \le \abs{\arg z} \le \pi \Longrightarrow \abs{E_{\alpha,1}(z)}\leq{C}({1+|z|})^{-1}.
\ee
Next, we pick  $\beta_0 \in \left( 0 , \min \left( \frac{\pi}{2} , \frac{\pi}{\alpha}-\frac{\pi}{2} \right) \right)$ and then choose $\beta \in (0, \pi)$ so small that
$   \abs{\arg(z-s)} \le \beta_0,\ z \in C_{T-\delta,\beta}=\{T-\delta +re^{i\theta}:\
r \in (0,\infty),\ \theta\in(-\beta,\beta)\},\ s \in (0,T-\delta).  $
Evidently, for all $z \in C_{T-\delta,\beta}$, all $s \in (0,T-\delta)$, and all $k \in \N$, we have
$ \abs{\arg(-\lambda_k (z-s)^\alpha)} \ge \pi - \alpha \beta_0 > \frac{\pi \alpha}{2}.$
So we get 
\bel{pp2b} 
\abs{E_{\alpha,1}(-\lambda_k(z-s)^\alpha)}\leq C\lambda_k^{-1} \abs{z-T+\delta}^{-\alpha},
\ee
upon applying \eqref{e2} with $\mu= \pi - \alpha \beta_0$. 
Consequently,
$$F(z)=\sum_{k=1}^\infty \frac{(\phi_{k}'(0))^2}{\lambda_k} \left(\int_0^{T-\delta}E_{\alpha,1}(-\lambda_{k}(z-s)^{\alpha})g'(s)\d s\right) $$
is well-defined for all $z\in \mathcal C_{T-\delta,\beta}$, by Lemma \ref{l1}. Moreover, according to \eqref{e3}, we have
\bel{e5}
\partial_xu(0,t)=F(t),\quad t\in(T-\delta,T).
\ee
It suffices to show that $F$ is holomorphic in $C_{T-\delta,\beta}$.  
By the holomorphicity of  $E_{\alpha,1}(z)$ in $\C$, the map 
$z\mapsto \int_0^{T-\delta}E_{\alpha,1}(-\lambda_{k}(z-s)^{\alpha})g'(s)\d s$, $k \in \N, $
is holomorphic in $\mathcal C_{T-\delta,\beta_\star}$. Therefore, for any compact subset $K$ of $\mathcal C_{T-\delta,\beta}$, it is enough to show that the sequence of functions
\bel{pp2c} F_N(z)= \sum_{k=1}^N \frac{(\phi_{k}'(0))^2}{\lambda_k} \left(\int_0^{T-\delta}E_{\alpha,1}(-\lambda_{k}(z-s)^{\alpha})g'(s)\d s\right),\quad N \in \mathbb N, \ee
converges uniformly in $K$. Indeed, by \eqref{pp2b}, we get for all $z\in K$ and all $M,N\in\mathbb{N}$ with $M< N$, that
\begin{align}
\abs{F_N(z)-F_M(z)} &\le \sum_{k=M+1}^N \frac{(\phi_{k}'(0))^2}{\lambda_k} \left(\int_0^{T-\delta}|E_{\alpha,1}(-\lambda_{k}(z-s)^{\alpha})||g'(s)|\d s\right) \nonumber\\
&\le C \left( \sum_{k=M+1}^N \frac{(\phi_{k}'(0))^2}{\lambda_k^2} \right) \norm{g'}_{L^1(0,T-\delta)} \left(\sup_{z\in K} \abs{z-T-\delta}^{-\alpha}\right) . \label{e4}
\end{align}
Note that we have $\sup_{z\in K} \abs{z-T-\delta}^{-\alpha}<\infty$ because $K$ is a compact subset of $\mathcal C_{T-\delta,\beta}$.
Thus, we deduce from \eqref{e4} and Lemma \ref{l1} that
$ \lim_{M,N\to+\infty} \sup_{z\in K} \abs{F_N(z)-F_M(z)}=0,$
showing that $(F_N)$ is a uniformly Cauchy sequence in $K$. Therefore, it converges uniformly in $K$ to $F$, which is holomorphic in $K$. Since $K$ is arbitrary in $\mathcal C_{T-\delta,\beta}$, $F$ is holomorphic in $\mathcal C_{T-\delta,\beta}$, and the desired assertion follows readily from this, \eqref{e5} and the inclusion $(T-\delta,\infty) \subset \mathcal C_{T-\delta,\beta}$.
\end{proof}

By repeating the above argumentation we can show a similar result for the solution $w$ of problem \eqref{eq2}. 

\begin{proposition}\label{pp3} Let the condition of Proposition \ref{pp1} and \eqref{q} be fulfilled. Then problem \eqref{eq2} admits a unique weak solution $w\in L^1(0,T;H^2(0,\ell))$. Moreover, if condition \eqref{t2b} is fulfilled, then the restriction of the map $t\mapsto\partial_xw(0,t)$ to $(T-\delta,T)$ admits an analytic extension to $(T-\delta,+\infty)$.
\end{proposition}

\section{Proof of Theorems \ref{t1} and \ref{t2}}\label{sec:proof}
Now we present the detailed proofs of Theorems \ref{t1} and \ref{t2}. The key tools in the analysis include solution representation and analyticity in Section \ref{sec:prelim}, Laplace transform and an improved version of Gel'fand-Levitan theory for the inverse Sturm-Liouville problem due to Sini \cite{Sini:2004}.

\subsection{Proof of Theorem \ref{t1}}

We divide the lengthy proof of the theorem into four steps.

\medskip

\noindent{\it Step 1: Analytic extension.}
For $j=1,2$, we denote by $A_j$ the operator $A$ defined in Section \ref{sec-prem} with $(\ell,\rho,q)=(\ell_j,\rho_j,q_j)$. That is, $A_j$ is the linear operator in $L^2(0,\ell_j;\rho_j \d x)$ acting as $\rho_j^{-1}(-\partial_x^2 u+q_j)$ on its domain $D(A_j)=H_0^1(0,\ell_j)\cap H^2(0,\ell_j)$.  We denote 
by $(\lambda_{j,k})_{k\in\mathbb N}$ the increasing sequence of simple eigenvalues of $A_j$, and pick an $L^2(0,\ell_j;\rho_j \d x)$-orthonormal basis $(\phi_{j,k})_{k \in \mathbb N}$ of eigenfunctions of $A_j$ such that
$ A_j \phi_{j,k} = \lambda_{j,k} \phi_{j,k},$ $ k \in \mathbb \N. $ 
Then by Proposition \ref{pp1},
we have
\bel{e9}
\partial_x u_j(0,t)= v_j^\prime(0) g(t) +\sum_{k=1}^\infty \frac{(\phi_{j,k}'(0))^2}{\lambda_{j,k}} \left(\int_0^t E_{\alpha_j,1}(-\lambda_{j,k}(t-s)^{\alpha_j})g'(s)\d s\right),\quad  t\in(0,T),
\ee
where $v_j$ is the $H^2(0,\ell_j)$-solution to \eqref{BVP} with $(\ell,q)=(\ell_j,q_j)$. 
Further, since $g$ vanishes in $(T-\delta,T)$ according to \eqref{t1a}, the above identity yields 
\bel{e10}
\partial_xu_j(0,t)=\sum_{k=1}^\infty \frac{(\phi_{j,k}'(0))^2}{\lambda_{j,k}} \left(\int_0^{T-\delta}  E_{\alpha_j,1}(-\lambda_{j,k}(t-s)^{\alpha_j})g'(s)\d s\right),\quad  t\in(T-\delta,T).
\ee
Now, by  Proposition \ref{pp2}, we can extend $t \mapsto \partial_x u_j(0,t)$ into an analytic function in $(T,\infty)$ and denote by $h_j$ the extension over $(0,\infty)$. In light of \eqref{e10}, we have
\bel{e11} 
h_j(t)=\sum_{k=1}^\infty \frac{(\phi_{j,k}'(0))^2}{\lambda_{j,k}} \left(\int_0^{T-\delta}E_{\alpha_j,1}(-\lambda_{j,k}(t-s)^{\alpha_j})g'(s)\d s\right),\quad t\in(T-\delta,\infty).
\ee
Next, since $t_n \in (T-\delta,T)$ for all $n \in \mathbb N$, we deduce from \eqref{t1b} and \eqref{e10}-\eqref{e11} that
$
h_1(t_n)=h_2(t_n)$, $ n\in\mathbb N.$ Since $(t_n)$ has an accumulation point at $T$ and $h=h_2-h_1$ is analytic in $(T-\delta,\infty)$, this entails 
$h(t)=0$ for  $t\in(T-\delta,+\infty),$
by the isolated zeros principle. Therefore, $h$ is supported on $[0,T-\delta]$ and its Laplace transform $\hat{h}$ is 
an entire function in $\mathbb C$. 

\medskip
\noindent {\it Step 2: Laplace transform.} 
Since $g \in W^{1,1}(0,T)$ is supported in $[0,T-\delta]$ and satisfies $g(0)=0$, it is extendable to a $W^{1,1}(0,\infty)$ function, still denoted by $g$, upon setting $g(s)=0$ for all $s \in (T,\infty)$. Then, in light of \eqref{e9} and \eqref{e11}, we find
$$ \begin{aligned}h_j(t)&=v_j^\prime(0) g(t) +\sum_{k=1}^\infty \frac{(\phi_{j,k}'(0))^2}{\lambda_{j,k}} \left(\int_0^t E_{\alpha_j,1}(-\lambda_{j,k}(t-s)^{\alpha_j})g'(s)\d s\right)\\
&=v_j^\prime(0) g(t) +\sum_{k=1}^\infty \frac{(\phi_{j,k}'(0))^2}{\lambda_{j,k}}\left(  E_{\alpha_j,1}(-\lambda_{j,k} s^{\alpha_j}) \ast g' \right) (t) ,\quad  t \in (0,\infty),\end{aligned} $$
where the symbol $\ast$ denotes the convolution operation on $[0,+\infty)$  between two functions. 
For all $k\in\mathbb N$, by setting 
$h_{j,k}(t)= \left(  E_{\alpha_j,1}(-\lambda_{j,k} s^{\alpha_j}) \ast g' \right) (t)$, for $t \in (0,\infty)$ and applying \cite[Eq. (1.80)]{P}, we have
$$\hat{h}_{j,k}(p)=\frac{p^{\alpha_j-1}\widehat{g'}(p)}{p^{\alpha_j}+\lambda_{j,k}}=\frac{p^{\alpha_j}\hat{g}(p)}{p^{\alpha_j}+\lambda_{j,k}},\quad p\in\mathbb C,\ \re (p)>\lambda_{j,k}^{\frac{1}{\alpha_j}},$$
where $\re(p)$ stands for the real part of $p$.
Since the maps
$p\mapsto \hat{h}_{j,k}(p)$ and $ p\mapsto \frac{p^{\alpha_j}\hat{g}(p)}{p^{\alpha_j}+\lambda_{j,k}}$
are holomorphic in $p\in\mathbb C_+:=\{z\in\mathbb C:\ \Re(z)>0\}$, we have
\bel{add1}\hat{h}_{j,k}(p)=\frac{p^{\alpha_j}\hat{g}(p)}{p^{\alpha_j}+\lambda_{j,k}},\quad p\in\mathbb C_+,\ k\in\mathbb N,\ j=1,2.\ee
In addition, from \eqref{essess} we obtain
$$\abs{\sum_{k=1}^N\frac{(\phi_{j,k}'(0))^2}{\lambda_{j,k}}h_{j,k}(t)}\leq C\sum_{k=1}^N\frac{(\phi_{j,k}'(0))^2}{\lambda_{j,k}^2}\int_0^t(t-s)^{-\alpha_j}|g'(s)|\d s,\quad t\in\R_+,\ N\in\mathbb N,$$
and, in view of Lemma \ref{l1} and \eqref{add1}, for all $p\in\mathbb C_+$, we find
\begin{align*}
    &\abs{e^{-pt}\sum_{k=1}^N\frac{(\phi_{j,k}'(0))^2}{\lambda_{j,k}}h_{j,k}(t)}\\
    \leq &C\left(\sum_{k=1}^{\infty}\frac{(\phi_{j,k}'(0))^2}{\lambda_{j,k}^2}\right)\int_0^te^{-\re(p)(t-s)}(t-s)^{-\alpha_j}e^{-\re(p) s}|g'(s)|\d s,\  t\in\R_+,\ N\in\mathbb N.
\end{align*}
Since $g'\in L^1(\R_+)$, by Young's inequality for convolution product and Lebesgue's dominated convergence theorem, we get
\bel{add2} 
\hat{h}_j(p)=v_j^\prime(0)\hat{g}(p) +\sum_{k=1}^\infty \frac{(\phi_{j,k}'(0))^2}{\lambda_{j,k}}\hat{h}_{j,k}(p)=v_j^\prime(0) \hat{g}(p)+\hat{g}(p)\sum_{k=1}^\infty\frac{(\phi_{j,k}'(0))^2p^{\alpha_j}}{\lambda_{j,k}(p^{\alpha_j}+\lambda_{j,k})},\quad p\in\mathbb C_+.
\ee
Put $\mu_1=\min((\frac{\lambda_{1,1}}{2})^{\frac{1}{\alpha_1}},(\frac{\lambda_{1,2}}{2})^{\frac{1}{\alpha_2}})>0$. Since $(\lambda_{j,k})_{k\in\mathbb N}$ is an increasing sequence, for all $k\in\mathbb N$ and $j=1,2$, we have
$\mu_1^{\alpha_j}\leq \frac{\lambda_{j,1}}{2}\leq\frac{\lambda_{j,k}}{2}$.
Next, setting $D_{\mu_1}:=\{z\in\mathbb C:\ |z|<\mu_1,\ z\notin(-\infty,0]\}$, we get for all $k\in\mathbb N$ and for $j=1,2$, 
$$\abs{\frac{(\phi_{j,k}'(0))^2}{\lambda_{j,k}(p^{\alpha_j}+\lambda_{j,k})}}\leq \frac{(\phi_{j,k}'(0))^2}{\lambda_{j,k}(\lambda_{j,k}-\mu_1^{\alpha_j})}\leq \frac{2(\phi_{j,k}'(0))^2}{\lambda_{j,k}^2},\quad p\in D_{\mu_1}. $$
Thus, by Lemma \ref{l1}, the sequence
$\sum_{k=1}^N\frac{(\phi_{j,k}'(0))^2}{\lambda_{j,k}(p^{\alpha_j}+\lambda_{j,k})}$, $ N\in\mathbb N$,
converges uniformly with respect to $p\in D_{\mu_1}$. It follows from this, \eqref{add2} and the holomorphicity of the  Laplace transform 
$\hat{g}$ in $\mathbb C$ (since $g$ is a compactly supported function), that $\hat{h}_j(p)$ can be extended holomorphically to $D_{\mu_1}$, as
\bel{e12} 
\hat{h}_j(p)=\hat{g}(p)v_j^\prime(0) +\hat{g}(p)\sum_{k=1}^\infty\frac{p^{\alpha_j}(\phi_{j,k}'(0))^2}{\lambda_{j,k}(p^{\alpha_j}+\lambda_{j,k})},\quad p\in D_{\mu_1}.
\ee
Further, by taking $p=R e^{\pm i \theta}$ in \eqref{e12}, with  fixed $R \in (0,\mu_1)$ and $\theta \in (0,\pi)$, and then sending $\theta$ to $\pi$, we get 
\begin{align}
 \lim_{\theta \to \pi} \left( \hat{h}_j(R e^{i \theta}) - \hat{h}_j(R e^{-i \theta}) \right)
= & 2i \sin(\alpha_j\pi) R^{\alpha_j} \hat{g}(-R)
\sum_{k=1}^\infty\frac{(\phi_{j,k}'(0))^2}{R^{2 \alpha_j}+ 2 R^{\alpha_j} \cos(\alpha_j \pi) \lambda_{j,k}+\lambda_{j,k}^2}.
\label{e13}
\end{align}
Therefore, as $R \downarrow 0$, we have
\begin{align}
 \lim_{\theta \to \pi} \left( \hat{h}_j(R e^{i \theta}) - \hat{h}_j(R e^{-i \theta}) \right) 
= &2i  \sin(\alpha_j\pi) R^{\alpha_j} \hat{g}(-R) \left( c_j - 2\cos(\alpha_j\pi) c_j' R^{\alpha_j} +{\mathcal O}(R^{2 \alpha_j}) \right), \label{e14}
\end{align}
where 
\bel{e15}
c_j=\sum_{k=1}^\infty \frac{(\phi_{j,k}'(0))^2}{\lambda_{j,k}^2}\quad \mbox{and}\quad c_j'=\sum_{k=1}^\infty \frac{(\phi_{j,k}'(0))^2}{\lambda_{j,k}^3}. 
\ee
Note that we have 
$c_j' \le \frac{c_j}{\lambda_{j,1}}$ and hence $c_j'<\infty$, from Lemma \ref{l1}.\\

\noindent {\it Step 3: Identification of the fractional order.}
Since $h=h_2-h_1$, we deduce from Step 2 that
$$\hat{h}(p)=\int_0^{+\infty}e^{-pt}(h_2(t)-h_1(t))dt=\hat{h}_2(p)-\hat{h}_1(p),\quad p\in\mathbb C_+,$$
where the map $\hat{h}_j(p)$, $j=1,2$, is given by \eqref{add2}. Combining this with \eqref{e12} and using the holomorphicity of $\hat{h}$ in $\mathbb C$ implies that this identity still holds on  $D_{\mu_1}$  and 
$$\hat{h}(p)=\hat{h}_2(p)-\hat{h}_1(p),\quad p\in D_{\mu_1},$$
where the map $\hat{h}_j(p)$ for $j=1,2$ is given by \eqref{e12}. Moreover, we have
$$ \lim_{\theta \to \pi} \left( \hat{h}(Re^{i\theta}) - \hat{h}(Re^{-i\theta}) \right) = \hat{h}(-R)-\hat{h}(-R) = 0,\quad R \in (0,\mu_1), $$
from the continuity of the (entire) function $\hat{h}$ at $p=-R$, and it follows that
\bel{e15b} \lim_{\theta \to \pi} \left( \hat{h}_2(R e^{i \theta}) - \hat{h}_2(R e^{-i \theta}) \right) - \lim_{\theta \to \pi} \left( \hat{h}_1(R e^{i \theta}) - \hat{h}_1(R e^{-i \theta}) \right) = 0,\quad R \in (0,\mu_1).
\ee
By combining this with \eqref{e14}, we get
\begin{align} 
& \sin(\alpha_1\pi) R^{\alpha_1} \hat{g}(-R) \left( c_1 - 2\cos(\alpha_1\pi) c_1' R^{\alpha_1} +{\mathcal O}(R^{2 \alpha_1}) \right)
 \nonumber \\
= &
\sin(\alpha_2\pi) R^{\alpha_2} \hat{g}(-R)
\left( c_2 - 2\cos(\alpha_2\pi) c_2' R^{\alpha_2} +{\mathcal O}(R^{2 \alpha_2}) \right),\quad \mbox{as }R \downarrow 0.
\label{e16}
\end{align}
Further, since $g$ is a non-zero compactly supported function, its Laplace transform 
$\hat{g}$ is holomorphic and not everywhere zero in $\C$. Thus, there exists a decreasing sequence $(r_n)_n$ of positive real numbers lying in $(0,\mu_1)$ and converging to $0$, such that
$\hat{g}(-r_n) \neq 0$ for all $ n\in\mathbb N.$
Thus, taking $R=r_n$ in \eqref{e16} gives 
\begin{align} 
& \sin(\alpha_1\pi) r_n^{\alpha_1} \left( c_1 - 2\cos(\alpha_1\pi) c_1' r_n^{\alpha_1} +{\mathcal O}(r_n^{2 \alpha_1}) \right)
 \nonumber \\
= &
\sin(\alpha_2\pi) r_n^{\alpha_2} 
\left( c_2 - 2\cos(\alpha_2\pi) c_2' r_n^{\alpha_2} +{\mathcal O}(r_n^{2 \alpha_2}) \right),\quad \mbox{as }n \to \infty.
\label{e17}
\end{align}
Let us assume  that $\alpha_1 < \alpha_2$. Then, multiplying both sides of \eqref{e17} by $r_n^{-\alpha_1}$ and sending $n$ to infinity yield
$\sin(\alpha_1\pi)c_1=0$, and hence $c_1=0$, since $\alpha_1 \pi \in (0,\pi)$. In light of \eqref{e15} and the fact that $\lambda_{1,k}>0$ for all $k \in \mathbb N$, this leads to the contradiction $\phi_{1,k}'(0)=0$. Thus we have 
$\alpha_1 \ge \alpha_2$ and hence $\alpha_1=\alpha_2$, as $\alpha_1$ and $\alpha_2$ play symmetric roles.\\

\noindent {\it Step 4: End of the proof.} In view of {\it Step 3}, we set $\alpha=\alpha_1=\alpha_2$ in \eqref{e15b} and then infer from 
\eqref{e13} that
$$ \hat{g}(-R)\sum_{k=1}^\infty \frac{(\phi_{1,k}'(0))^2}{(R^{\alpha}+\lambda_{1,k}e^{i\alpha\pi})(R^{\alpha}+\lambda_{1,k}e^{-i\alpha\pi})} = \hat{g}(-R)\sum_{k=1}^\infty \frac{(\phi_{2,k}'(0))^2}{(R^{\alpha}+\lambda_{2,k}e^{i\alpha\pi})(R^{\alpha}+\lambda_{2,k}e^{-i\alpha\pi})},\quad R \in (0,\mu_1). 
$$
Next, since $\hat{g}$ is holomorphic and not everywhere zero in $\C$, there exist two positive real numbers $r_1<r_2<\mu_1$ such that $\hat{g}(-R)\neq0$ for all $R\in(r_1,r_2)$. Thus, the above identity implies
\bel{t1g}
\sum_{k=1}^\infty\frac{(\phi_{1,k}'(0))^2}{(R^{\alpha}+\lambda_{1,k}e^{i\alpha\pi})(R^{\alpha}+\lambda_{1,k} e^{-i\alpha\pi})}=\sum_{k=1}^\infty \frac{(\phi_{2,k}'(0))^2}{(R^{\alpha}+\lambda_{2,k}e^{i\alpha\pi})(R^{\alpha}+\lambda_{2,k}e^{-i\alpha\pi})},\quad R\in(r_1,r_2).
\ee
Meanwhile, in light of Lemma \ref{l1}, one can deduce that the map
$$z \mapsto\sum_{k=1}^\infty\frac{(\phi_{j,k}'(0))^2}{(z+\lambda_{j,k}e^{i\alpha\pi})(z+\lambda_{j,k}e^{-i\alpha\pi})}$$
is meromorphic in $\mathbb C$ with simple poles at $\lambda_{j,k}e^{i (1\pm \alpha)\pi}$, $k\in\mathbb N$. Then, \eqref{t1g} implies
\bel{t1h}
\sum_{k=1}^\infty\frac{(\phi_{1,k}'(0))^2}{(ze^{i\alpha\pi}+\lambda_{1,k})(ze^{-i\alpha\pi}+\lambda_{1,k})}=\sum_{k=1}^\infty\frac{(\phi_{2,k}'(0))^2}{(ze^{i\alpha\pi}+\lambda_{2,k})(ze^{-i\alpha\pi}+\lambda_{2,k})},\ee
for all $z \in \mathbb C\setminus \{\lambda_{j,k}e^{i(1\pm \alpha)\pi},\ j=1,2,\ k\in\mathbb N\}$. This yields 
$(\phi_{1,k}'(0))^2=(\phi_{2,k}'(0))^2\ \mbox{and}\ \lambda_{1,k}=\lambda_{2,k}$ for $ k\in\mathbb N,$
or equivalently 
$
(|\phi_{1,k}'(0)|,\lambda_{1,k})=(|\phi_{2,k}'(0)|,\lambda_{2,k})$ for $ k\in\mathbb N.$
It follows from this and \cite[Theorem 1]{Sini:2004}, that
$(\ell_1,\rho_1,q_1)=(\ell_2,\rho_2,q_2),$
which yields the desired result.

\begin{remark} Note that one  beneficial memory effect of fractional diffusion equations lies in the discontinuities on the whole negative real axis $(-\infty,0)$ of the Laplace transform in time of solutions of \eqref{eq1}, which is key to the proof of Theorem \ref{t1}. This property, which is due to the fact that $\alpha\not\in\mathbb N$, allows us to derive \eqref{e13} and to obtain formulas \eqref{e16} and \eqref{t1g} from \eqref{e15b}. Clearly this property is no longer valid when $\alpha_1=\alpha_2=1$ since in that case we would instead have
$$\lim_{\theta \to \pi} \hat{h}_j(R e^{i \theta}) - \hat{h}_j(R e^{-i \theta})  = 0,\quad R \in (0,\mu_1),\ j=1,2.$$\end{remark}

\subsection{Proof of Theorem \ref{t2}}
For $j=1,2$, consider the operator $B_j=\rho_j^{-1}(-\partial_x^2 u+q_j)$ acting on $L^2(0,\ell_j;\rho_j\d x)$ with its domain $D(B_j)=\{ v\in H^2(0,\ell_j):\ v(0)=v'(\ell_j)=0\}$. The spectrum of $B_j$ consists of an increasing sequence $(\mu_{j,k})_{k\in\mathbb N}$ of simple and positive eigenvalues. We associate with these eigenvalues an orthonormal basis $(\psi_{j,k})_{k\in\mathbb N}$ of $L^2(0,\ell_j;\rho_j\d x)$. The solution $w$ of \eqref{eq2}, with $\alpha=\alpha_j$, $\ell=\ell_j$, $q=q_j$ and $\rho=\rho_j$, takes the form
$$w_j(x,t)=-\sum_{k=1}^\infty \left(\int_0^t(t-s)^{\alpha_j-1}E_{\alpha_j,\alpha_j}(-\mu_{j,k}(t-s)^{\alpha_j})h(s)\d s\right)\psi_{j,k}'(0)\psi_{j,k}(x),\quad x\in(0,\ell_j),\ t\in(0,T).$$
Using this representation and mimicking the proof of Theorem \ref{t2} yield $\alpha_1=\alpha_2$ and
$|\psi_{1,k}'(0)|=|\psi_{2,k}'(0)|$, $\mu_{1,k}=\mu_{2,k}$ for all $k\in\mathbb N.$
Then, by combining \cite[Section 5]{Sini:2004} with \cite[Theorem 1]{Sini:2004}, we deduce \eqref{t2b}.

\section{Numerical results and discussions}
\label{sec:numer}
In this section, we present numerical results to demonstrate the feasibility of the simultaneous recovery of unknown parameters and the interval size.
\subsection{Levenberg-Marquardt method}
First we describe the numerical algorithm for obtaining the reconstruction. Since the inverse problem involves four unknown parameters ($\alpha$, $\ell$, $\rho$ and $q$), it is numerically very challenging to recover all of them, and also the recovery of the order $\alpha$ has been extensively studied (see, e.g., \cite{JinKian:2021,LiaoWei:2019}). Thus, we focus on recovering the potential $q$ and the interval length $\ell$. In Examples \ref{ex:Dirichlet} and \ref{ex:Neumann}, we assume that the order $\alpha$ and the density $\rho$ are given, whereas in Example \ref{ex:Dirichlet_unknown_rho}, we assume that only the order  $\alpha$ is given. Numerically, since these parameters have different influence on the measurement, standard iterative methods, e.g. Landweber method and conjugate gradient, do not work very well. In several studies (see, e.g., \cite{LiaoWei:2019,JinKianZhou:2023}, the Levenberg-Marquardt method \cite{Levenberg:1944,Marquardt:1963} has shown excellent numerical performance for solving related inverse problems. Thus, we employ it for the numerical reconstruction below. 

We describe the method for recovering $q$ and $\ell$, and the extension to more general cases is direct. We define a nonlinear operator $F:(q,\ell)\in L^2(0,\ell)\times \mathbb{R}_+\rightarrow \partial_x u(0,t)\in L^2(T_0, T)$, where the function $u$ solves problem \eqref{eq1} or \eqref{eq2} with parameters $q$ and $\ell$. Fix an initial guess $(q^0,\ell^0)$. Now given $(q^k,\ell^k)$, we find the next approximation $(q^{k+1},\ell^{k+1})$ by 
\begin{equation*}
    (q^{k+1},\ell^{k+1})=\arg\min J_k(q,\ell),
\end{equation*}
with the functional $J_k(v,\ell)$ (based at $(q^k,\ell^k)$) given by
\begin{align*}
    J_k(q,\ell)=&\|F(q^{k},\ell^k)-z^{\delta}+\partial_qF(q^k,\ell^k)(q-q^k)+\partial_{\ell} F(q^k,\ell^k)(\ell-\ell^k) \|^2_{L^2(T_0,T)}\\
    &+\beta_q^k\|q-q^k\|_{L^2(0,\ell^k)}^2+{\beta_{\ell}^k}|\ell-\ell^k|^2 +{\mu}\|q \|_{L^2(0,\ell^k)}^2 ,
\end{align*}
where $z^\delta$ is the noisy data, $[T_0,T]$ is the measurement time horizon, $\beta_q^k$, $\beta_\ell^k$ and $\mu$ are positive scalars, and $\partial_q F (q^k,\ell^k)$ and $\partial_\ell F (q^k,\ell^k)$ are the Jacobians of the forward map $F$ in $q$ and $\ell$, respectively. The parameters $\beta_q^k$ and $\beta_\ell^k$ are often decreased geometrically with decreasing factors $\gamma_q,\gamma_\ell\in(0,1)$: $\beta_q^{k+1}=\gamma_q\beta_q^k $ and $\beta_\ell^{k+1}=\gamma_\ell \beta_\ell^k $. The parameter $\mu$ is fixed during the iteration. The derivative $\partial_q F(q^k,\ell^k)$  can be evaluated explicitly. Indeed, let  $u^k$ solve problem \eqref{eq1} with  $(q,\ell)=(q^k,\ell^k)$. Then the directional derivative is $\partial_qF(q^k,\ell^k)[h]=\partial_x w(0,t) $, with
\begin{equation*} 
    \left\{ 
\begin{aligned}
&\rho \partial_{t}^{\alpha}w-\partial_{x}^2 w+q^k w    =-h u^k,  &&\mbox{in }(0,\ell^k)\times (0,T],\\
&w(\ell^k,t)=0,\, w(0,t)=0, && t\in (0,T)\\
&w=0, &&\mbox{in }(0,\ell^k)\times \{0\}.
\end{aligned}
\right.
\end{equation*}
Furthermore, $\partial_{\ell}F(q^k,\ell^k)$ can be computed using the difference quotient $\partial_{\ell}F(q,\ell)\approx(\delta \ell)^{-1}(F(q,\ell+\delta \ell)-F(q,\ell))$, where $\delta\ell$ is a small number, fixed at $\delta\ell=\text{1e-3}$ below. 
 
\subsection{Numerical results and discussions}
Now we present numerical  results by the Levenberg-Marquardt method. Throughout, we take the domain $\Omega=(0,1)$, i.e. $\ell^\dag=1$, and $T=1$. We discretize  problems \eqref{eq1} and \eqref{eq2} using the Galerkin method with conforming piecewise linear finite element in space and backward Euler convolution quadrature in time \cite{JinZhou:2023book}. In all experiments, we take a mesh size $h=1/100$ and time step size $\tau=1/1000$. The accuracy of a recovered potential $q$ is measured by the relative $L^2(0,\ell)$ error $e_q=\|q-q^\dag\|_{L^2(0,\ell )}/\|q^\dag\|_{L^2(0,\ell)}$.
The residual $r$ of the recovered tuple $(q,\ell)$  is $r(q,\ell)=\|F(q,\ell)-z^\delta\|_{L^2(T_0,T)}$. The exact data $\partial_x u^\dag(0,t)$ is generated on a fine space-time mesh, and the noisy $z^\delta$ is generated as 
\begin{equation*}
    z^\delta(t)=\partial_x u^\dag(0,t)+\epsilon \| \partial_x u^\dag(0,\cdot)\|_{L^\infty(T_0,T)} \xi (t),
\end{equation*}
where $\xi(t)$ follows the standard Gaussian noise and $\epsilon>0$ denotes the relative noise level.

\begin{example}\label{ex:Dirichlet}
This example is for the inverse problem of equation \eqref{eq1}. Let the boundary data $g(t)=\chi_{[0,0.5]}$, where $\chi$ denotes the characteristic function. Fix $\alpha$, and $\rho(x)=1+0.5\chi_{[0.5,\ell]}$. The potential $q^{\dag}(x)=10x(1-x)^2$ and $\ell^{\dag}=1$ are unknown. The measurement data $z^\delta$ is collected in interval $[T_0,T]=[0.6,1]$.
\end{example}

\begin{figure}[hbt!]
    \centering
    \setlength{\tabcolsep}{0pt}
    \begin{tabular}{ccc}
        \includegraphics[width=0.33\textwidth]{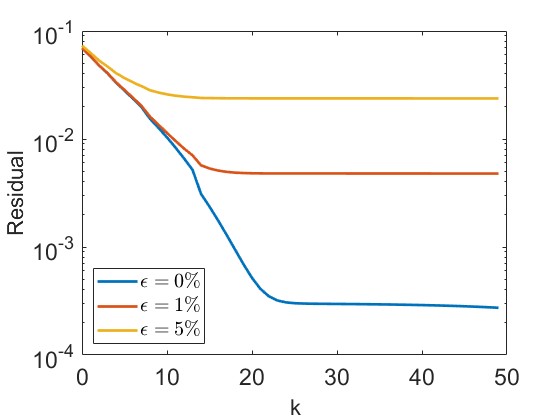} &
        \includegraphics[width=0.33\textwidth]{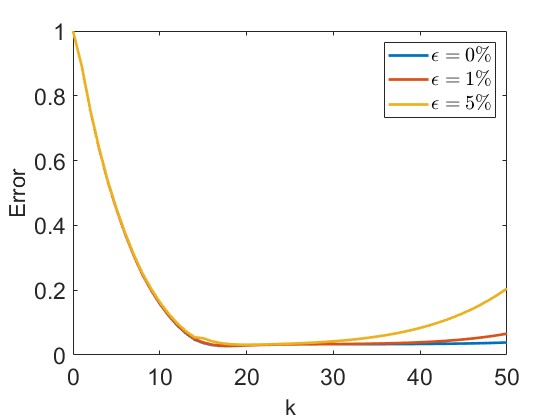} &
        \includegraphics[width=0.33\textwidth]{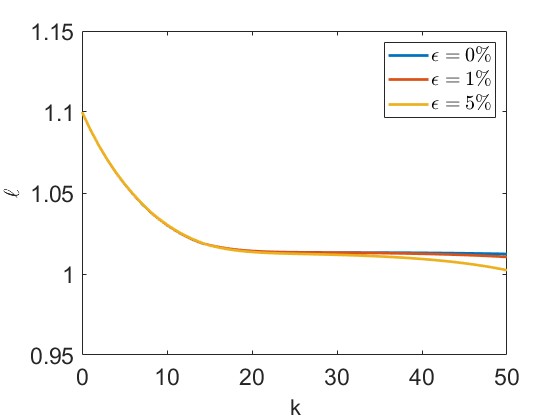} \\
        (a) residual $r$ & (b) error $e_q$ & (c) interval length $\ell$
    \end{tabular}
    \caption{The convergence of the algorithm for Example \ref{ex:Dirichlet},  at three noise levels, $\alpha=0.75$.}
\label{fig:Dirichlet_convergence}
\end{figure}

\begin{figure}[hbt!]
    \centering
    \setlength{\tabcolsep}{0pt}
    \begin{tabular}{ccc}
    \includegraphics[width=0.33\textwidth]{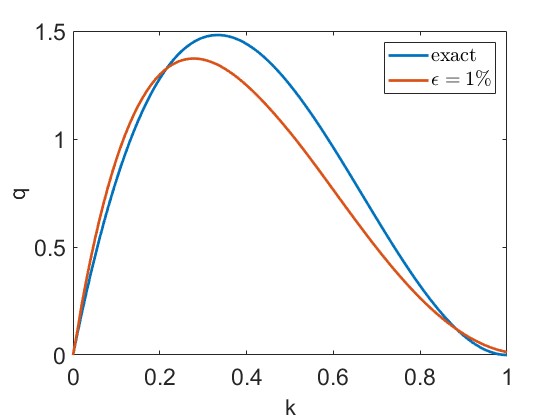} &
    \includegraphics[width=0.33\textwidth]{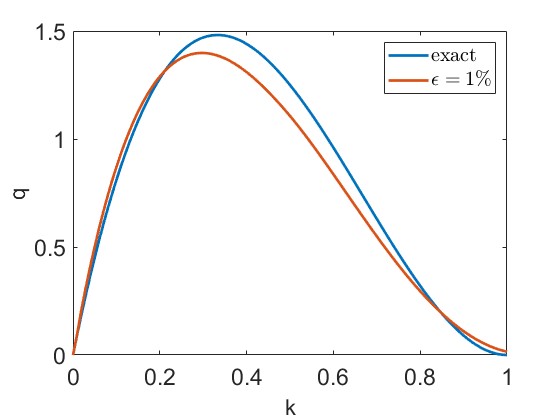} &
    \includegraphics[width=0.33\textwidth]{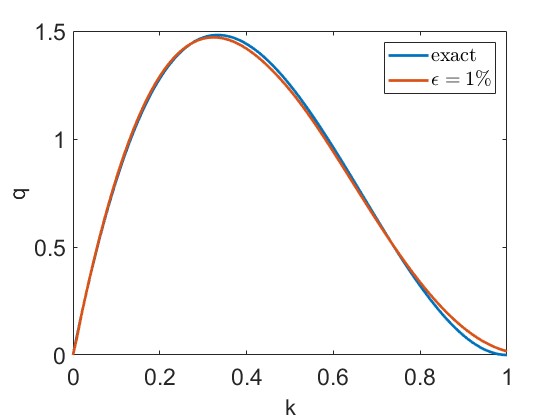} \\
        (a) $\alpha=0.25$ & (b) $\alpha=0.50$ & (c) $\alpha=0.75$
    \end{tabular}
    \caption{The reconstructions of the potential $q$ for Example \ref{ex:Dirichlet}.}
    \label{fig:Dirichlet_reconstruction}
\end{figure}

\begin{figure}[hbt!]
    \centering
    \setlength{\tabcolsep}{0pt}
    \begin{tabular}{ccc}
        \includegraphics[width=0.33\textwidth]{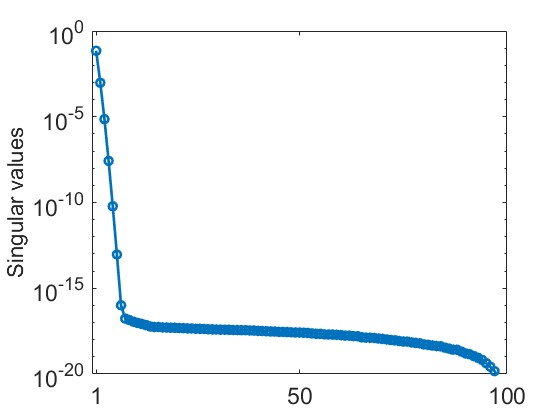} &
        \includegraphics[width=0.33\textwidth]{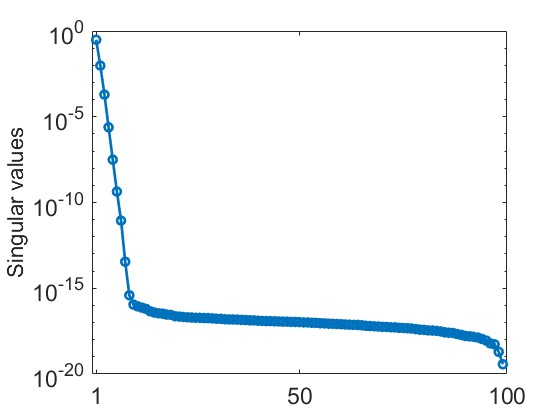} &
        \includegraphics[width=0.33\textwidth]{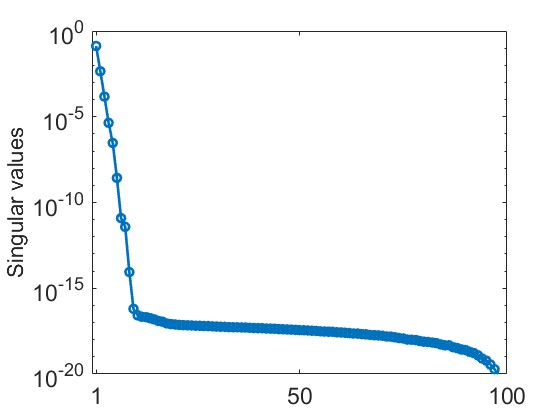} 
    \end{tabular}
    \caption{The singular values of the Jacobian $\partial_q F$ in  Examples \ref{ex:Dirichlet}, \ref{ex:Neumann}, and \ref{ex:Dirichlet_unknown_rho} (from left to right).}
    \label{fig:Singularvalues}
\end{figure}

In Fig. \ref{fig:Dirichlet_convergence}, we show the convergence of the Levenberg-Marquardt method. We initialize the method at $q^0(x)=0$ and $\ell^0=1.1$. For exact data $\epsilon=0\%$, the residual $r$ decays to $\text{2e-4}$, due to the presence of discretization error, and the error $e_q$ decays to $\text{3e-2}$ and then stabilizes. For noisy data, we observe a typical semi-convergence behavior: The error $e_q$ first decreases, levels off and eventually starts to increase (after about twenty iterations). This behavior is also observed for the interval length $\ell$. Thus, early stopping is critical to get reasonable reconstructions. The semi-convergence phenomenon is more pronounced at the higher noise level  $\epsilon=5\%$. Since the Levenberg-Marquardt parameters $\beta_q$ and $\beta_\ell$ decays geometrically, the reconstruction eventually blows up as the number of iterations increases, due to the ill-posedness of the inverse problem. Indeed, Fig. \ref{fig:Singularvalues}(a) shows that the singular values of the linearized operator $\partial_q F$ decay extremely fast to $10^{-15}$ (the machine precision). This indicates that the inverse problem is severely ill-posed and one can only expect to capture a few singular functions (corresponding to large singular values).
Fig. \ref{fig:Dirichlet_reconstruction} shows the recovered potential $q$ with $1\%$ noise data, where the results are chosen such that relative error $e_q$ is smallest along the iteration trajectory; see also Table \ref{table:Dirichlet} for quantitative results. Due to the severe ill-posedness of the inverse problem, we have taken relatively large parameters $\beta_q$ and $\beta_\ell$ such that the initial iterations are relatively stable.

\begin{table}[hbt!]
  \centering
  \begin{threeparttable}
    \caption{The parameter choice and numerical results for Example \ref{ex:Dirichlet} at noise level $1\%$.\label{table:Dirichlet}}
    \begin{tabular}{cccccccc} 
        \toprule
            $\alpha$ & $\beta_q^0$ & $\beta_\ell^0$ & $\gamma_q$ & $\gamma_\ell$ & $\mu$ & $e_q$  & $\ell$ \\
        \midrule
            $0.25$ & 2e-2 & $200$ & $0.9$ & $0.9$ & 1e-9 & 1.12e-1 & $1.0133$\\
            $0.50$ & 5e-2 & $500$ & $0.9$ & $0.9$ & 1e-9 & 7.63e-2  & $1.0129$\\
            $0.75$ & 1e-1 & $1000$ & $0.9$ & $0.9$ & 1e-9 & 2.81e-2  & $1.0151$\\
        \bottomrule
    \end{tabular}
      \end{threeparttable}
\end{table}

\begin{example}\label{ex:Neumann}
This example is for the inverse problem of equation \eqref{eq2}. Let the boundary data $h(t)=\chi_{[0,0.8]}$. Fix $\alpha$, and $\rho(x)\equiv 1$. The potential $q^{\dag}(x)=\frac{1}{1+e^{-10x}}$ and $\ell^{\dag}=1$ are unknown. The measurement data $z^\delta$ is collected in interval $[T_0,T]=[0.9,1]$.
\end{example}

\begin{figure}[hbt!]
    \centering
    \setlength{\tabcolsep}{0pt}
    \begin{tabular}{ccc}
        \includegraphics[width=0.33\textwidth]{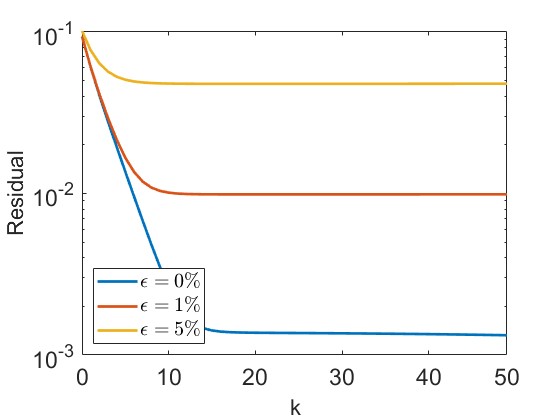} &
        \includegraphics[width=0.33\textwidth]{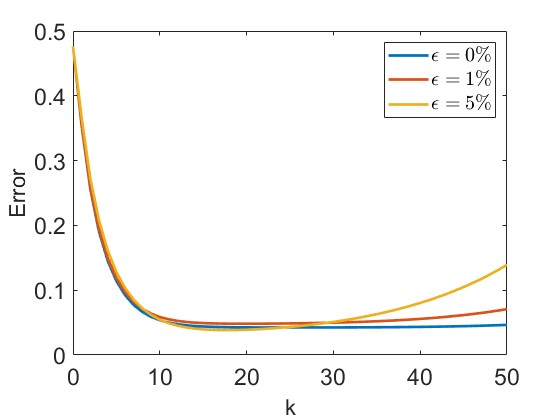} &
        \includegraphics[width=0.33\textwidth]{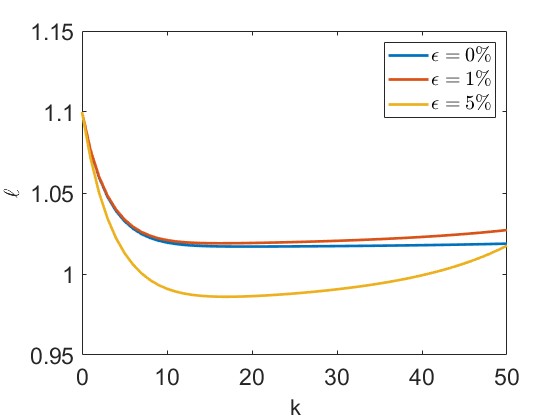} \\
        (a) residual $r$ & (b) error $e_q$ & (c) interval length $\ell$
    \end{tabular}
    \caption{The convergence of the algorithm for Example \ref{ex:Neumann} at three noise levels, for $\alpha=0.75$.}
    \label{fig:Neumann_convergence}
\end{figure}

\begin{figure}[hbt!]
    \centering
    \setlength{\tabcolsep}{0pt}
    \begin{tabular}{ccc}
        \includegraphics[width=0.33\textwidth]{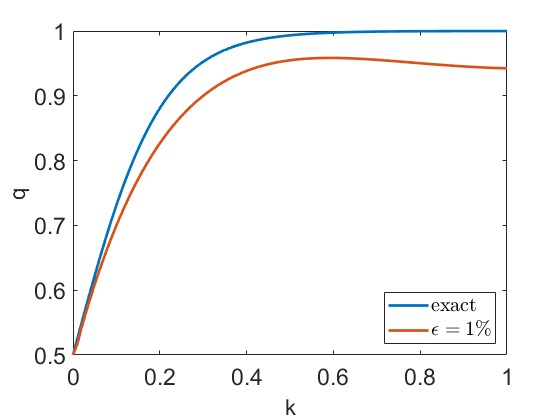} &
        \includegraphics[width=0.33\textwidth]{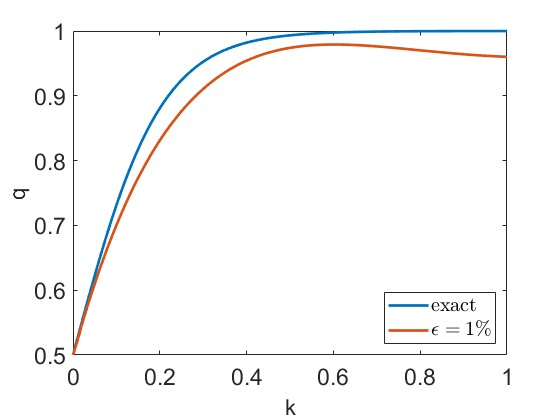} &
        \includegraphics[width=0.33\textwidth]{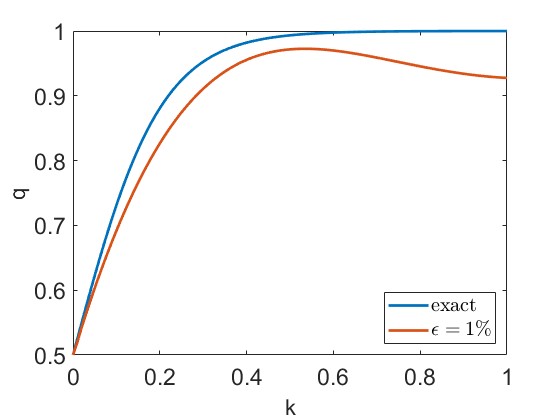} \\
        (a) $\alpha=0.25$ & (b) $\alpha=0.50$ & (c) $\alpha=0.75$
    \end{tabular}
    \caption{The reconstructions of the potential $q$ for Example \ref{ex:Neumann} with $1\%$ noise.}
    \label{fig:Neumann_reconstruction}
\end{figure}

The numerical results for Example \ref{ex:Neumann} are shown in Fig. \ref{fig:Neumann_convergence}. The initial guess is taken as $q^0(x)=0.5$ and $\ell^0=1.1$. The convergence plot is quite similar to Example \ref{ex:Dirichlet}: there exists a  semi-convergence phenomenon of the iterates, and the error $e_q$ first decreases and then increases afterwards. The numerical reconstructions of the potential $q$ are given in Fig. \ref{fig:Neumann_reconstruction}; see Table \ref{table:Neumann} for quantitative results. The reconstruction quality does not depend much on the fractional order $\alpha$, which is consistent with the theoretical results. 

\begin{table}[hbt!]
  \centering
  \begin{threeparttable}
  \caption{The parameter choice and numerical results for Example \ref{ex:Neumann} at noise level $1\%$.\label{table:Neumann}}
    \begin{tabular}{cccccccc} 
        \toprule
            $\alpha$ & $\beta_q^0$ & $\beta_\ell^0$ & $\gamma_q$ & $\gamma_\ell$ & $\mu$ & $e_q$  & $\ell$ \\
        \midrule
            $0.25$ & 2e-2 & $20$ & $0.9$ & $0.95$ & 1e-6 & 4.83e-2 & $0.9998$\\
            $0.50$ & 5e-2 & $40$ & $0.9$ & $0.95$ & 1e-6 & 3.23e-2  & $0.9826$\\
            $0.75$ & 1e-1 & $80$ & $0.9$ & $0.95$ & 1e-6 & 4.84e-2  & $1.0191$\\
        \bottomrule
    \end{tabular}
    \end{threeparttable}
\end{table}

The last example aims at recovering also a constant $\rho$.
\begin{example}\label{ex:Dirichlet_unknown_rho}
This example is for the inverse problem of equation \eqref{eq1}. Fix $\alpha$. The density $\rho^{\dag}(x)= 1$, the potential $q^{\dag}(x)=10x(1-x)^2$ and $\ell^{\dag}=1$ are unknown.
\end{example}

\begin{figure}[hbt!]
    \centering
    \setlength{\tabcolsep}{0pt}
    \begin{tabular}{cc}
        \includegraphics[width=0.45\textwidth]{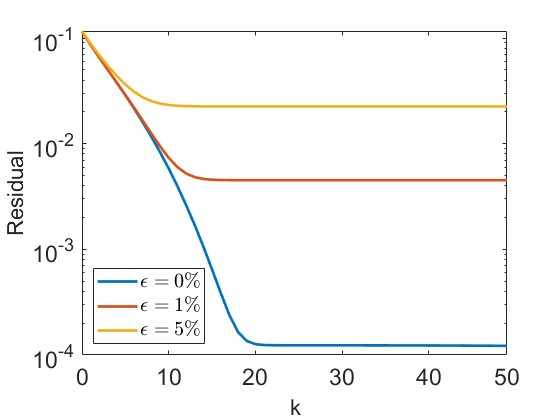} &
        \includegraphics[width=0.45\textwidth]{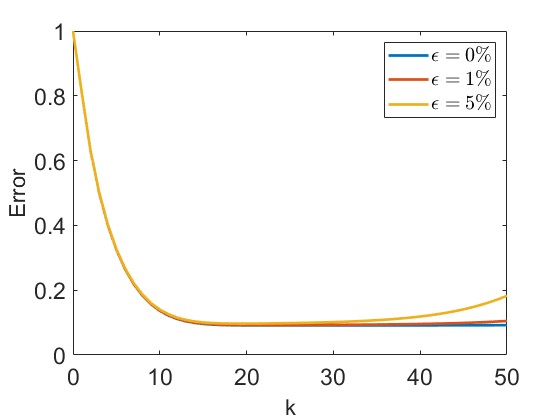} \\
        (a) residual $r$ & (b) error $e_q$ \\
        \includegraphics[width=0.45\textwidth]{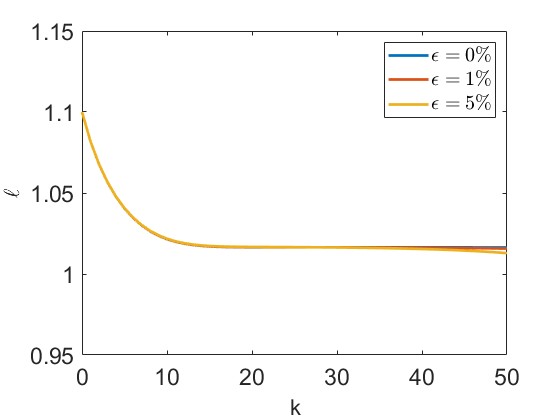} & 
        \includegraphics[width=0.45\textwidth]{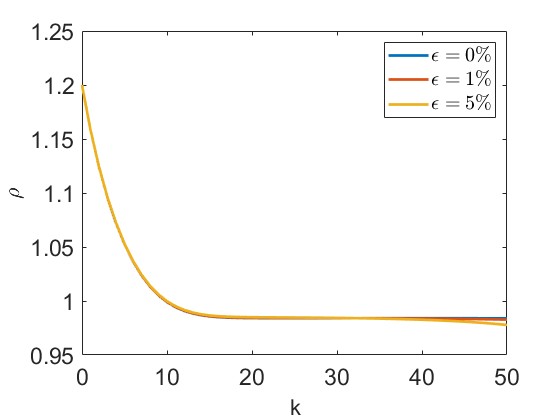} \\
        (c) interval length $\ell$ & (d) $\rho$ 
    \end{tabular}
    \caption{The convergence of the algorithm for Example \ref{ex:Dirichlet_unknown_rho} at three noise levels, for $\alpha=0.75$.}
    \label{fig:Dirichlet_unknown_rho_convergence}
\end{figure}

\begin{figure}[hbt!]
    \centering
        \setlength{\tabcolsep}{0pt}
    \begin{tabular}{ccc}
        \includegraphics[width=0.33\textwidth]{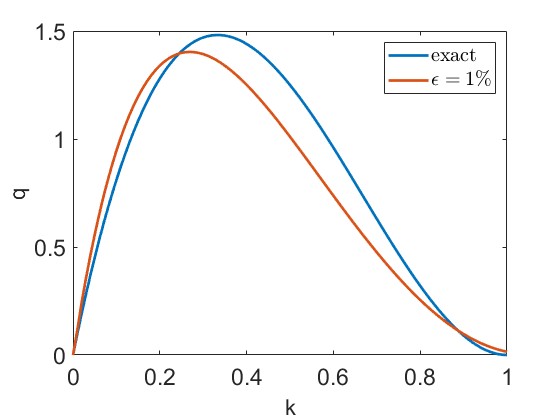} &
        \includegraphics[width=0.33\textwidth]{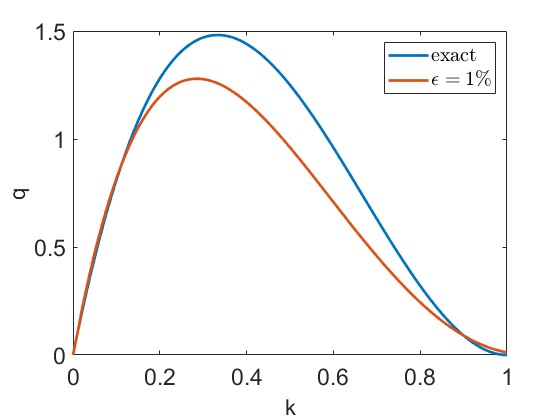} &
        \includegraphics[width=0.33\textwidth]{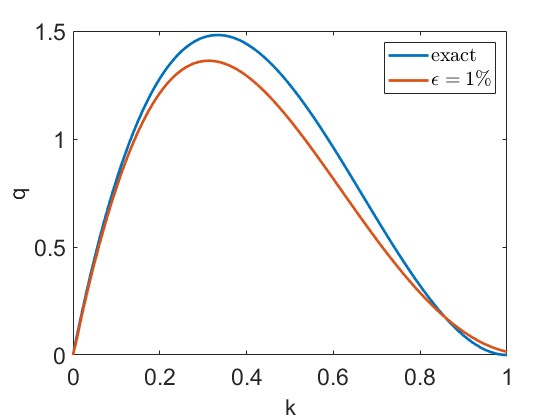} \\
        (a) $\alpha=0.25$ & (b) $\alpha=0.50$ & (c) $\alpha=0.75$
    \end{tabular}
    \caption{The reconstructions of the potential $q$ for Example \ref{ex:Dirichlet_unknown_rho} with $1\%$ noise .}
    \label{fig:Dirichlet_unknown_rho_reconstruction}
\end{figure}

The simultaneous recovery of the density $\rho$ along with $q$ and $\ell$ is very challenging. So we focus on the case that $\rho$ is constant. In Fig. \ref{fig:Dirichlet_unknown_rho_convergence}, we show the convergence of the algorithm with the initial guess $q^0(x)=0$, $\ell^0=1.1$ and $\rho^0=1.2$. Like before, the method exhibits a semi-convergence behavior, especially for noisy data. Since the density $\rho$ is constant, the convergence of $\rho$ appears to be more robust to the iteration number $k$ than the potential $q$. The choice of parameters and quantitative numerical results are given in Table \ref{table:Dirichlet_unknown_rho} and  Fig. \ref{fig:Dirichlet_unknown_rho_reconstruction}. Compared with the results for Example \ref{ex:Dirichlet}, the relative error $e_q$ is larger due to the unknown density $\rho$, indicating the challenges of  recovering multiple coefficients.

\begin{table}[hbt!]
  \centering
  \begin{threeparttable}
  \caption{The parameter choice and numerical results for Example \ref{ex:Dirichlet_unknown_rho} with noise level $\epsilon=1\%$.\label{table:Dirichlet_unknown_rho}}
    \begin{tabular}{ccccccccccc} 
        \toprule
            $\alpha$ & $\beta_q^0$ & $\beta_\ell^0$ & $\beta_\rho^0$ & $\gamma_q$ & $\gamma_\ell$ &  $\gamma_\rho$ & $\mu$ & $e_q$  & $\ell$ & $\rho$ \\
        \midrule
            $0.25$ & 2e-2 & $200$ & $50$ & $0.9$ & $0.9$ & $0.9$ & 1e-9 & 1.39e-1 & $1.0061$ & $1.0081$\\
            $0.50$ & 5e-2 & $500$ & $100$ & $0.9$ & $0.9$ & $0.9$ & 1e-9 & 1.66e-1 & $1.0169$ & $0.9865$\\          
            $0.75$ & 1e-1 & $1000$ & $200$ & $0.9$ & $0.9$ & $0.9$ & 1e-9 & 9.24e-2 & $1.0165$ & $0.9843$\\
        \bottomrule
    \end{tabular}
  \end{threeparttable}
\end{table}

These results show clearly the feasibility of recovering multiple parameters in the 1D subdiffusion model from the lateral flux measurement, confirming the theoretical predictions in Theorems \ref{t1} and \ref{t2}.

\section*{Acknowledgments}
The work of B. Jin is supported by Hong Kong RGC General Research Fund (Project
14306423), and a start-up fund from The Chinese University of Hong Kong. The work of Z. Zhou is supported by Hong Kong
Research Grants Council (15303122) and an internal grant of Hong Kong Polytechnic University (Project ID: P0038888, Work
Programme: 1-ZVX3).
Rachid Zarouf acknowledges the Theoretical Physics Center (CPT) from Aix-Marseille University for the kind hosting during his collaboration with the other authors of the article. 

\pagebreak
\bibliographystyle{abbrv}
\bibliography{frac}

\begin{thebibliography}{10}

\bibitem{AdamsGelhar:1992}
E.~E. Adams and L.~W. Gelhar.
\newblock Field study of dispersion in a heterogeneous aquifer: 2. spatial
  moments analysis.
\newblock {\em Water Res. Research}, 28(12):3293--3307, 1992.

\bibitem{AS}
S.~Avdonin and T.~I. Seidman.
\newblock Identification of {{\(q(x)\)}} in {{\(u_ t=\Delta u-qu\)}} from
  boundary observations.
\newblock {\em SIAM J. Control Optim.}, 33(4):1247--1255, 1995.

\bibitem{CenJinZhou:2023}
S.~Cen, B.~Jin, Y.~Liu, and Z.~Zhou.
\newblock Recovery of multiple parameters in subdiffusion from one lateral
  boundary measurement.
\newblock {\em Inverse Problems}, 39(10):104001, 31, 2023.

\bibitem{CY}
J.~Cheng and M.~Yamamoto.
\newblock Identification of convection term in a parabolic equation with a
  single measurement.
\newblock {\em Nonlinear Anal., Theory Methods Appl., Ser. A, Theory Methods},
  50(2):163--171, 2002.

\bibitem{HatanoHatano:1998}
Y.~Hatano and N.~Hatano.
\newblock Dispersive transport of ions in column experiments: An explanation of
  long-tailed profiles.
\newblock {\em Water Res. Research}, 34(5):1027--1033, 1998.

\bibitem{HLLZ}
T.~Helin, M.~Lassas, L.~Ylinen, and Z.~Zhang.
\newblock Inverse problems for heat equation and space-time fractional
  diffusion equation with one measurement.
\newblock {\em J. Differ. Equations}, 269(9):7498--7528, 2020.

\bibitem{JK}
J.~Janno and Y.~Kian.
\newblock Inverse source problem with a posteriori boundary measurement for
  fractional diffusion equations.
\newblock {\em Math. Methods Appl. Sci.}, 46(14):15868--15882, 2023.

\bibitem{Jin:2021}
B.~Jin.
\newblock {\em {Fractional Differential Equations}}.
\newblock Springer-Nature, Switzerland, 2021.

\bibitem{JinKian:2021}
B.~Jin and Y.~Kian.
\newblock Recovery of the order of derivation for fractional diffusion
  equations in an unknown medium.
\newblock {\em SIAM J. Appl. Math.}, 82(3):1045--1067, 2022.

\bibitem{JinKianZhou:2023}
B.~Jin, Y.~Kian, and Z.~Zhou.
\newblock Inverse problems for subdiffusion from observation at an unknown
  terminal time.
\newblock {\em SIAM J. Appl. Math.}, 83(4):1496--1517, 2023.

\bibitem{JinRundell:2015}
B.~Jin and W.~Rundell.
\newblock A tutorial on inverse problems for anomalous diffusion processes.
\newblock {\em Inverse Problems}, 31(3):035003, 40, 2015.

\bibitem{JinZhou:2021ip}
B.~Jin and Z.~Zhou.
\newblock Recovering the potential and order in one-dimensional time-fractional
  diffusion with unknown initial condition and source.
\newblock {\em Inverse Problems}, 37(10):105009, 28, 2021.

\bibitem{JinZhou:2023book}
B.~Jin and Z.~Zhou.
\newblock {\em {Numerical Treatment and Analysis of Time-Fractional Evolution
  Equations}}, volume 214 of {\em Applied Mathematical Sciences}.
\newblock Springer, Cham, 2023.

\bibitem{JingYamamoto:2023}
X.~Jing and M.~Yamamoto.
\newblock Simultaneous uniqueness for multiple parameters identification in a
  fractional diffusion-wave equation.
\newblock {\em Inverse Probl. Imaging}, 16(5):1199--1217, 2022.

\bibitem{Kian:2022ip}
Y.~Kian.
\newblock Simultaneous determination of different class of parameters for a
  diffusion equation from a single measurement.
\newblock {\em Inverse Problems}, 38(7):075008, 29, 2022.

\bibitem{KianLiLiuYamamoto:2021}
Y.~Kian, Z.~Li, Y.~Liu, and M.~Yamamoto.
\newblock The uniqueness of inverse problems for a fractional equation with a
  single measurement.
\newblock {\em Math. Ann.}, 380(3-4):1465--1495, 2021.

\bibitem{KY1}
Y.~Kian and M.~Yamamoto.
\newblock Well-posedness for weak and strong solutions of non-homogeneous
  initial boundary value problems for fractional diffusion equations.
\newblock {\em Fract. Calc. Appl. Anal.}, 24(1):168--201, 2021.

\bibitem{KilbasSrivastavaTrujillo:2006}
A.~A. Kilbas, H.~M. Srivastava, and J.~J. Trujillo.
\newblock {\em Theory and {A}pplications of {F}ractional {D}ifferential
  {E}quations}.
\newblock Elsevier Science B.V., Amsterdam, 2006.

\bibitem{KJ}
N.~Kinash and J.~Janno.
\newblock An inverse problem for a generalized fractional derivative with an
  application in reconstruction of time- and space-dependent sources in
  fractional diffusion and wave equations.
\newblock {\em Mathematics}, 2019(7):1138, 2019.

\bibitem{Kou:2008}
S.~C. Kou.
\newblock Stochastic modeling in nanoscale biophysics: subdiffusion within
  proteins.
\newblock {\em Ann. Appl. Stat.}, 2(2):501--535, 2008.

\bibitem{KubicaYamamoto:2020}
A.~Kubica, K.~Ryszewska, and M.~Yamamoto.
\newblock {\em {Time-Fractional Differential Equations---a Theoretical
  Introduction}}.
\newblock Springer, Singapore, 2020.

\bibitem{Levenberg:1944}
K.~Levenberg.
\newblock A method for the solution of certain non-linear problems in least
  squares.
\newblock {\em Quart. Appl. Math.}, 2:164--168, 1944.

\bibitem{LiYamamoto:2019coeff}
Z.~Li and M.~Yamamoto.
\newblock Inverse problems of determining coefficients of the fractional
  partial differential equations.
\newblock In {\em {Handbook of Fractional Calculus with Applications. {V}ol.
  2}}, pages 443--464. De Gruyter, Berlin, 2019.

\bibitem{LiaoWei:2019}
K.~Liao and T.~Wei.
\newblock Identifying a fractional order and a space source term in a
  time-fractional diffusion-wave equation simultaneously.
\newblock {\em Inverse Problems}, 35(11):115002, 23, 2019.

\bibitem{Marquardt:1963}
D.~W. Marquardt.
\newblock An algorithm for least-squares estimation of nonlinear parameters.
\newblock {\em J. Soc. Indust. Appl. Math.}, 11:431--441, 1963.

\bibitem{Meerschaert:2019}
M.~M. Meerschaert and A.~Sikorskii.
\newblock {\em {Stochastic Models for Fractional Calculus}}, volume~43.
\newblock De Gruyter, Berlin, second edition, 2019.

\bibitem{MetzlerJeon:2014}
R.~Metzler, J.-H. Jeon, A.~G. Cherstvy, and E.~Barkai.
\newblock Anomalous diffusion models and their properties: non-stationarity,
  non-ergodicity, and ageing at the centenary of single particle tracking.
\newblock {\em Phys. Chem. Chem. Phys.}, 16(44):24128--24164, 2014.

\bibitem{MetzlerKlafter:1998}
R.~Metzler, J.~Klafter, and I.~M. Sokolov.
\newblock Anomalous transport in external fields: continuous time random walks
  and fractional diffusion equations extended.
\newblock {\em Phys. Rev. E}, 58(2):1621--1633, 1998.

\bibitem{Nigmatulin:1986}
R.~R. Nigmatullin.
\newblock The realization of the generalized transfer equation in a medium with
  fractal geometry.
\newblock {\em Phys. Stat. Sol. B}, 133:425--430, 1986.

\bibitem{P}
I.~Podlubny.
\newblock {\em Fractional differential equations}.
\newblock Academic Press, Inc., San Diego, CA, 1999.

\bibitem{PT}
J.~P{\"o}schel and E.~Trubowitz.
\newblock {\em Inverse {S}pectral {T}heory}, volume 130.
\newblock Academic Press, New York, 1987.

\bibitem{RundellYamamoto:2018}
W.~Rundell and M.~Yamamoto.
\newblock Recovery of a potential in a fractional diffusion equation.
\newblock Preprint, arXiv:1811.05971, 2018.

\bibitem{RundellYamamoto:2023}
W.~Rundell and M.~Yamamoto.
\newblock Uniqueness for an inverse coefficient problem for a one-dimensional
  time-fractional diffusion equation with non-zero boundary conditions.
\newblock {\em Appl. Anal.}, 102(3):815--829, 2023.

\bibitem{SY}
K.~Sakamoto and M.~Yamamoto.
\newblock Initial value/boundary value problems for fractional diffusion-wave
  equations and applications to some inverse problems.
\newblock {\em J. Math. Anal. Appl.}, 382(1):426--447, 2011.

\bibitem{Sini:2004}
M.~Sini.
\newblock On the one-dimensional {G}elfand and {B}org-{L}evinson spectral
  problems for discontinuous coefficients.
\newblock {\em Inverse Problems}, 20(5):1371--1386, 2004.

\bibitem{Yamamoto:2023isp}
M.~Yamamoto.
\newblock Uniqueness for inverse source problems for fractional diffusion-wave
  equations by data during not acting time.
\newblock {\em Inverse Problems}, 39(2):024004, 20, 2023.

\end{thebibliography}
\end{document}